\documentclass[journal,twoside,web]{ieeecolor}
\usepackage{generic}
\usepackage{cite}
\usepackage{amsmath,amssymb,amsfonts}
\usepackage{algorithmic}
\usepackage{graphicx}
\usepackage{algorithm,algorithmic}
\usepackage{epstopdf}

\newtheorem{assumption}{Assumption}
\newtheorem{theorem}{Theorem}

\allowdisplaybreaks
\usepackage{multirow}
\usepackage{subcaption}
\usepackage{lipsum}

\usepackage{textcomp}
\def\BibTeX{{\rm B\kern-.05em{\sc i\kern-.025em b}\kern-.08em
    T\kern-.1667em\lower.7ex\hbox{E}\kern-.125emX}}
\begin{document}

\title{Unbiased Extremum Seeking for MPPT in Photovoltaic Systems}

\author{Cemal Tugrul Yilmaz, Eric Foss, Mamadou Diagne, and Miroslav Krstic
\thanks{The authors are with the Department of Mechanical and Aerospace Engineering, University of California San Diego, La Jolla, CA 92093-0411 USA (e-mail: cyilmaz@ucsd.edu, efoss@ucsd.edu, mdiagne@ucsd.edu, krstic@ucsd.edu).}
}

\maketitle

\allowdisplaybreaks

\begin{abstract}
This paper presents novel extremum seeking (ES) strategies for maximum power point tracking (MPPT) in photovoltaic (PV) systems that ensure unbiased convergence and prescribed-time performance. Conventional ES methods suffer from steady-state bias due to persistent dither signal. We introduce two novel ES algorithms: the exponential unbiased ES (uES), which guarantees exponential convergence to the maximum power point (MPP) without steady-state oscillation bias, and the unbiased prescribed-time ES (uPT-ES), which ensures convergence within a user-defined time horizon. Both methods leverage time-varying perturbation amplitudes and demodulation gains, with uPT-ES additionally utilizing chirp signals to enhance excitation over finite-time intervals. Experimental results on a hardware-in-the-loop testbed validate the proposed algorithms, demonstrating improved convergence speed and tracking accuracy compared to classical ES, under both static and time-varying environmental conditions.

\end{abstract}

\begin{IEEEkeywords}
extremum seeking, prescribed-time control, maximum power point tracking, photovoltaic systems
\end{IEEEkeywords}

\section{Introduction}
\label{sec:introduction}

Energy today faces a multidimensional crisis involving rising costs, widespread energy poverty, serious environmental degradation, and escalating geopolitical risks such as oil-related conflicts and nuclear proliferation \cite{holdren1991population}. Managing the increasing energy demand driven by population growth and environmental concerns requires a balanced energy mix, emphasizing renewable sources like wind, hydroelectric, and solar power. According to the National Laboratory of Renewable Energy, in 2024, renewables accounted for $33\%$ of U.S. electric capacity and $24\%$ of generation. All non-carbon sources--including solar, wind, nuclear, hydropower, and geothermal--represent $41\%$ of energy capacity  and $40\%$ of energy  generation. Methodologies focused on optimizing renewable energy sources (RES), with particular emphasis on  photovoltaic distributed generation play a crucial role in overcoming a swiftly rising global energy demand.

\subsection{Photovoltaic MPPT Background}
Photovoltaic (PV) arrays exhibit highly nonlinear current–voltage (I-V) and power–voltage (P-V)
characteristics that depend strongly on irradiance and temperature.
For each irradiance and temperature condition, the P-V curve possesses a unique maximum power point (MPP).
Identifying this point in practice is a challenging task and requires maximum power point tracking (MPPT) controllers that continuously adjust the converter’s operating voltage to maximize power extraction \cite{kruitwagen2021global,verma2016maximum,ashwini2025,pandey2008high}. Conventional electronic MPPT algorithms like perturb-and-observe (PO) or incremental conductance are popular due to their simplicity and model-free operation \cite{ishaque2014performance, sera2013perturb}. A recent study proposed a robust modeling and control framework for solar photovoltaic conversion systems under real-time weather conditions \cite{sedraoui2025}. By modifying the PO algorithm and integrating a fixed-order $H_{\infty}$ controller, the method achieved superior tracking and robustness compared to standard PO-based MPPT, particularly under rapid meteorological changes. Another study proposed an adaptive MPPT framework based on model reference adaptive control (MRAC) to improve efficiency and robustness in photovoltaic systems \cite{manna2023design}. The MRAC-MPPT design outperformed conventional incremental conductance and perturb \& observe methods in tracking speed, accuracy, and stability under varying environmental conditions.

\subsection{Extremum Seeking in MPPT and Its Limitations}
As an advanced alternative to conventional MPPT methods, extremum seeking (ES) has been explored for maximum power point tracking in PV systems. ES is a dynamic, model-free optimization technique that injects small perturbation signals into the system to estimate the gradient of the function, guiding the system toward the optimum. Several ES-based MPPT schemes have been proposed in the literature. One of the earliest applications of ES to PV systems is presented in \cite{leyva2006mppt}. A notable development is the implementation of Newton-based ES in \cite{ghaffari2014power}, which enables user-defined convergence rates to the MPP, independent of system parameters, for cascaded PV architectures. This design is extended in \cite{li2014newton}, where the periodic probing signals are replaced with stochastic counterparts and tested for MPPT performance. In \cite{brunton2010maximum}, a different ES algorithm leverages the inverter’s natural ripple to track the MPP under rapidly varying irradiance conditions. To minimize steady-state perturbations while maintaining responsiveness to changing conditions, a Lyapunov-based switched ES approach is proposed and validated in \cite{moura2013lyapunov}. An adaptive amplitude ES scheme is introduced in \cite{tchouani2023climatic} to improve tracking performance under partial shading. Similarly, an adaptive ESC method has been employed for maximum power point tracking in photovoltaic systems, where the duty ratio of the DC–DC converter is optimized and radial basis function neural networks are used to approximate the nonlinear current–voltage characteristics \cite{li2011maximum}. An ESC-based MPPT method is developed in \cite{heydari2013new} to ensure global MPP tracking under partial shading, reducing ripple and improving efficiency without added complexity.  Additionally, a sliding-mode ES method is tested in alternator-based energy conversion systems to enhance robustness against external disturbances \cite{toloue2017multivariable}. A comprehensive overview and survey of ES-based MPPT strategies can be found in \cite{bizon2016global}. 

\begin{figure*}[t]
    \centering
    \includegraphics[width=1.8\columnwidth]{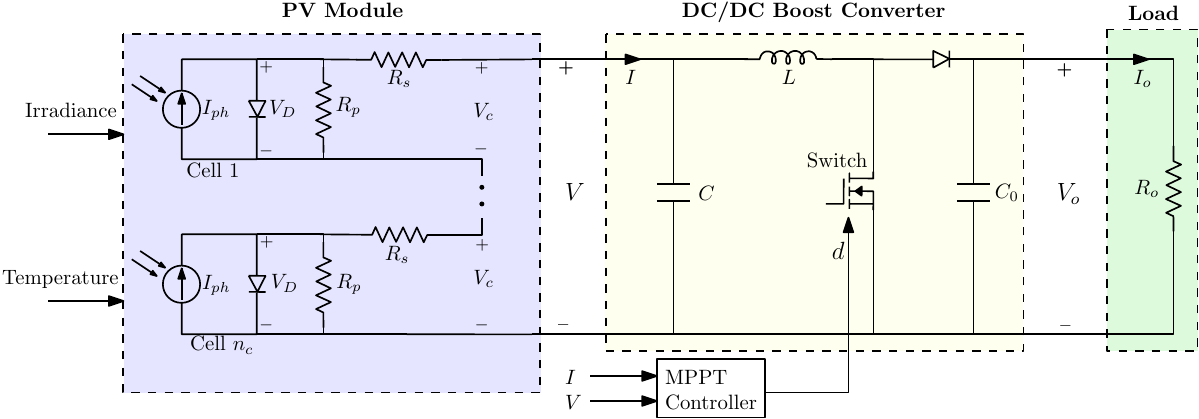}
    \caption{Block diagram of a PV system comprising a PV module, a DC/DC boost converter, and a load. The duty cycle $d$ is perturbed by an MPPT controller to maximize the output power $P = VI$, which is influenced by variations in solar irradiance and temperature.}
    \label{PVmoduleblock}
\end{figure*}

\subsection{Unbiased and Prescribed-Time Extremum Seeking for MPPT}
Efforts to reduce steady-state oscillations in ES algorithms have not only been pursued within the context of MPPT but have also attracted broader attention in the general ES literature. 
Although several designs achieve asymptotic convergence to the optimum without steady-state bias \cite{abdelgalil2021lie, guay2020uncertainty, haring2016asymptotic, suttner2019extremum}, the method proposed in \cite{YILMAZ2025112392} is the first to ensure exponential convergence to the optimum at a user-specified rate, and is referred to as unbiased extremum seeking (uES).
The key idea behind uES is to simultaneously reduce the perturbation amplitude, minimizing  exploration, while increasing the feedback gain to enhance exploitation of the gradient estimate. This trade-off drives the steady-state bias to zero exponentially. Building upon uES, prescribed-time uES (PT-uES) in \cite{YILMAZ2025112392} introduces time-varying gains and frequencies (e.g., chirp signals) to ensure convergence to the optimum within a user-defined finite time.

This paper makes three key contributions: $(i)$ the introduction of an exponentially converging unbiased ES for MPPT, $(ii)$ a novel prescribed-time variant that ensures convergence within a user-defined horizon, and $(iii)$ experimental validation demonstrating superior performance under static and dynamic irradiance.
To the best of our knowledge, no existing MPPT algorithm provides exponential, unbiased convergence to the MPP, nor does any offer a provable finite-time convergence guarantee. 
Despite the theoretical guarantee of unbiased convergence, practical implementations experience minor bias due to unavoidable noise and environmental variability.

The paper is organized as follows. Section II presents the PV system model and discusses the effects of irradiance and temperature. Section III outlines the problem statement. Section IV introduces the exponential unbiased extremum seeking algorithm for MPPT, and Section V extends this approach to the prescribed-time setting. Section VI details the experimental setup and results. Finally, Section VII concludes the paper.

\section{PV System Modeling}

In this section, we consider the modeling of a PV system consisting of a PV module, a DC/DC boost converter, and a load, as illustrated in Fig. \ref{PVmoduleblock}. We also discuss the effects of irradiance and temperature on the power–voltage characteristics of the PV cells, as well as the need for an MPPT algorithm to maximize the produced power.

\subsection{PV Module Model}

A PV module consists of multiple PV cells connected in series. The electrical behavior of a PV cell is commonly modeled using the single-diode representation, which captures the nonlinear relationship between the current $I$ and voltage $V$ as described by \begin{equation} 
I = I_{\mathrm{ph}} - I_0 \left( \exp\left( \frac{V_c + I R_s}{n V_T} \right) - 1 \right) - \frac{V_c + I R_s}{R_{p}}, \label{eq:PV_general} 
\end{equation} 
where the photogenerated current $I_{\mathrm{ph}}$ is approximated by 
\begin{equation} I_{\mathrm{ph}} = \left( I_{\mathrm{SC,ref}} + \alpha_I (T - T_{\mathrm{ref}}) \right) \frac{G}{G_{\mathrm{ref}}}, \label{eq:Iph} 
\end{equation} 
and the diode reverse saturation current $I_0$ is modeled as \begin{equation} I_0 = I_{0,\mathrm{ref}} \left( \frac{T}{T_{\mathrm{ref}}} \right)^3 \exp\left( \frac{q E_g}{n k_B} \left( \frac{1}{T_{\mathrm{ref}}} - \frac{1}{T} \right) \right). \label{eq:I0} 
\end{equation} 
Here, $R_s$ and $R_{p}$ denote the series and parallel resistances, respectively; $n$ is the diode ideality factor; $I_{\mathrm{SC,ref}}$ is the short-circuit current at reference irradiance $G_{\mathrm{ref}}$ and reference temperature $T_{\mathrm{ref}}$; $I_{\mathrm{o,ref}}$ is the reference saturation current, $\alpha_I$ is the temperature coefficient of the short-circuit current; $G$ is the incident solar irradiance; and $E_g$ is the semiconductor bandgap energy. The thermal voltage $V_T$ is given by $V_T = \frac{k_B T}{q}$, where $k_B$ is Boltzmann’s constant, $q$ is the elementary charge, and $T$ is the absolute temperature in Kelvin.

\subsection{DC-DC Boost Converter Model}

A DC-DC boost converter is commonly used in PV systems to connect the PV source to the load. It steps up the voltage and regulates the operating voltage/current through its duty cycle, which is controlled by an MPPT algorithm. In continuous conduction mode, where the inductor current remains positive throughout the switching cycle, the averaged continuous-time model of the boost converter is given by:
\begin{align}
    L \frac{dI_L}{dt} &= V - (1 - d) V_o, \label{eq:boost_L} \\
    C \frac{dV_o}{dt} &= (1 - d) I_L - I_0, \label{eq:boost_C}
\end{align}
where $L$ and $C$ are the inductance and capacitance,
$I_L$ is the inductor current, $V$ is the input voltage from the PV source, $V_o$ is the output voltage, $I_o$ is the output (load) current, and $d \in [0,1]$ is the duty cycle. At steady state, the time derivatives vanish, and \eqref{eq:boost_L} and \eqref{eq:boost_C} yield
\begin{align}
    V &= (1 - d) V_0, \label{eq:v_relation} \\
    I_{0} &= (1 - d)I. \label{eq:current_relation}
\end{align}

\subsection{Effect of Irradiance and Temperature}

\begin{figure}[t]
    \centering
    \begin{subfigure}[b]{\linewidth}
        \centering
        \includegraphics[width=\linewidth]{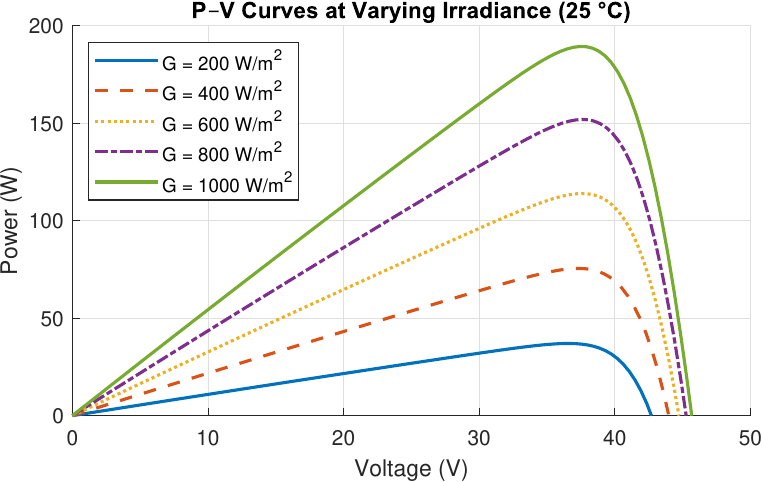}
        \caption{}
    \end{subfigure} \\ \vspace{0.2cm}
    \begin{subfigure}[b]{\linewidth}
        \centering
        \includegraphics[width=\linewidth]{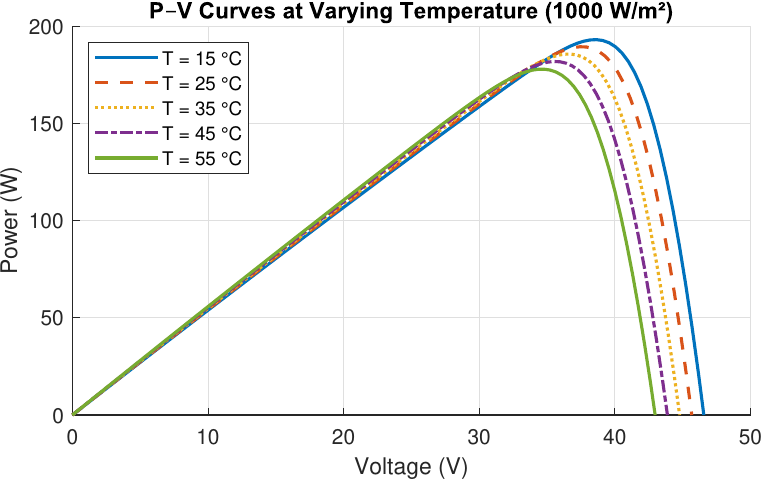}
        \caption{}
    \end{subfigure}
    \caption{P-V curves of a PV panel under (a) varying irradiance at 25$^\circ C$ and (b) varying temperature at 1000~$W/m^2$. Increased irradiance raises both power and optimal voltage, while higher temperature lowers the maximum power and shifts the optimal voltage leftward.}
    \label{fig:exponent_uES_track}
\end{figure}

The simulated power-voltage curves in Fig. \ref{fig:exponent_uES_track} show the performance of a silicon-based PV module under varying irradiance ($G$) and temperature ($T$). Parameters include $\alpha_{\mathrm{sc}} = 0.0047$ A/$^\circ$C, $a_{\mathrm{ref}} = 1.2$, $I_{\mathrm{L,ref}} = 5.5$ A, $I_{\mathrm{o,ref}} = 10^{-10}$ A, $R_{\mathrm{sh,ref}} = 200,\Omega$, $R_{\mathrm{s}} = 0.5,\Omega$, $E_{\mathrm{g,ref}} = 1.121$ eV, with $N_s = 60$ cells, $k = 1.38 \times 10^{-23} J/K$ , and $q = 1.60 \times 10^{-19}  C$.

The maximum power point lies between two extremes: the short-circuit condition, where the current is maximum and the voltage is zero in \eqref{eq:PV_general}, and the open-circuit condition, where the voltage is maximum and the current is zero in \eqref{eq:PV_general}. The short-circuit current $I_{\mathrm{sc}}$ and open-circuit voltage $V_{\mathrm{oc}}$ are commonly used to estimate the MPP in practice. The short-circuit current (approximately $I_{ph}$) increases nearly linearly with irradiance, while the open-circuit voltage has a logarithmic dependence on current and decreases with temperature. The plots show that increasing irradiance at $25^\circ$C boosts both $I_{\mathrm{sc}}$ and MPP, while increasing temperature at 1000 $W/m^2$ reduces $V_{\mathrm{oc}}$ and MPP, despite a slight rise in current.

\section{Problem Statement}

We consider the problem of regulating the duty cycle $ d \in [0,1] $ of a DC-DC boost converter to maximize the output power extracted from the PV panel. The power delivered to the load is given by  
\begin{equation}  
P = I V,  
\end{equation}  
where $I$ and $V$ denote the current and voltage at the PV terminals, respectively. PV modules exhibit nonlinear characteristics, as illustrated in Fig. \ref{fig:exponent_uES_track}. The power–voltage curve typically has a single peak corresponding to the maximum power point, which indicates the existence and uniqueness of an optimal operating voltage (and thus duty cycle) that maximizes output power.

Building on this observation and to enable comprehensive analysis and algorithm development, we introduce the following assumption on the map $P(d)$:
\begin{assumption} \label{assconvex}
The function $P$ is $\mathcal{C}^4$, and there exists $d^{*} \in \mathbb{R}$ such that
\begin{align}
    \frac{\partial}{\partial d} P(d^{*})={}&0, \label{derivativeP}\\
    \frac{\partial^2}{\partial d^2} P(d^{*})={}&h<0, \quad h=h^T. \label{doublederivativeP}
\end{align}
\end{assumption}
Assumption~\ref{assconvex} ensures that the function $P(d)$ has a unique local maximum at $ d^* $. The first-order condition~\eqref{derivativeP} implies stationarity, while the second-order condition~\eqref{doublederivativeP} guarantees strong concavity at the maximizer.

\section{Exponential Unbiased Extremum Seeker} \label{expExSect}

\begin{figure}[t]
    \centering
    \includegraphics[width=.9\columnwidth]{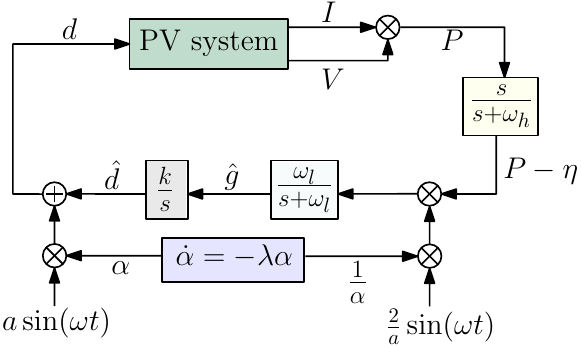} 
    \caption{Exponential uES scheme. The design employs an exponentially decaying function $\alpha$ to gradually diminish the effect of the perturbation signal $a\sin(\omega t)$, while its multiplicative inverse $\frac{1}{\alpha}$ correspondingly amplifies the effect of the demodulation signal $\frac{2}{a}\sin(\omega t)$.}
    \label{ESBlock}
\end{figure}

To seek the maximum power point through duty cycle regulation, we employ the exponential uES scheme, schematically illustrated in Fig.~\ref{ESBlock} and originally introduced in~\cite{YILMAZ2025112392}. The design parameters include the gain $k$, the high-pass and low-pass filter coefficients $\omega_h$ and $ \omega_l $, all of which are strictly positive real numbers. The perturbation amplitude $a$ is a real-valued constant.

The exponentially decaying function $ \alpha(t) $ is defined by
\begin{align}
    \dot{\alpha}(t) = -\lambda \alpha(t), \qquad \alpha(t_0) = \alpha_0, \label{mu}
\end{align}
where $ \lambda > 0 $ and $ \alpha_0 > 0 $ are design parameters.

To ensure convergence and stability of the closed-loop system, the following conditions must be satisfied:
\begin{align}
    \lambda &< \min \left\{ \frac{\omega_l}{2}, \frac{\omega_h}{2} \right\}, \label{cond1} \\
    k &> (\omega_l - \lambda) \frac{\lambda}{\omega_l} \left( \frac{1}{-h} \right) > 0, \label{cond2}
\end{align}

These conditions require that the perturbation amplitude $\alpha$ decays and the demodulation gain $1/\alpha$ increases at a controlled rate. To ensure convergence, the estimator $\hat{d}$ must adapt sufficiently quickly relative to this decay; in other words, the adaptation must outpace the loss of excitation.

We summarize closed-loop system depicted in Fig. \ref{ESBlock} as follows
\begin{align}
    \frac{d}{d t}\begin{bmatrix}
    \tilde{d} \\ \hat{g} \\ \tilde{\eta} \\ \alpha
    \end{bmatrix}=\begin{bmatrix}
        k \hat{g} \\ -\omega_l \hat{g}+\omega_l (P-P(d^*)-\tilde{\eta})\frac{1}{\alpha} \frac{2}{a} \sin(\omega t) \\
        -\omega_h \tilde{\eta}+\omega_h (P-P(d^*)) \\
        -\lambda \alpha
    \end{bmatrix}, \label{ffed}
\end{align}
in view of the transformations
\begin{align}
    \tilde{d}={}&\hat{d}-d^*, \label{transf1} \\
    \tilde{\eta}={}&\eta-P(d^*), \label{transf2}
\end{align}
where $\eta$ is governed by
\begin{align}
    \dot{\eta}=-\omega_h \eta+\omega_h P. \label{etadef}
\end{align}

The convergence result follows from \cite{YILMAZ2025112392}.
\begin{theorem} \label{the1}
Consider the feedback system \eqref{ffed} with the parameters 
that satisfy \eqref{cond1}, \eqref{cond2} under Assumption \ref{assconvex}. There exists $\bar{\omega}$ and for any $\omega > \bar{\omega}$ there exists an open ball $\mathcal{B}$ centered at the point $(\hat{d}, \hat{g}, {\eta}, \alpha)=(d^*, 0, P(d^*),0) = : \Upsilon$ such that for any initial condition starting in the ball $\mathcal{B}$, the system \eqref{ffed} has a unique solution and the solution converges exponentially to $\Upsilon$. Furthermore, $P(d(t))$ exponentially converges to $P(d^*)$.
\end{theorem}

A brief sketch of the proof is provided below:

\begin{proof} Let us proceed through the proof step by step.

\textbf{Step 1: State transformation.}
We introduce the following state transformations:
\begin{align}
     \tilde{d}_f=&{} \frac{1}{\alpha} \tilde{d}, \qquad
     \hat{g}_f={} \frac{1}{\alpha} \hat{g},  \qquad 
     \tilde{\eta}_f={} \frac{1}{\alpha^2} \tilde{\eta}, \label{trans}
\end{align}
which transform the system in \eqref{ffed} into:
\begin{align}
    &\frac{d}{d t}\begin{bmatrix}
    \tilde{d}_f & \hat{g}_f & \tilde{\eta}_f & \alpha
    \end{bmatrix}^T \nonumber \\
    &=\begin{bmatrix}
        \lambda \tilde{d}_f+k \hat{g}_f \\ (\lambda-{\omega}_l)\hat{g}_f +{\omega}_l\left(\nu(z)-\tilde{\eta}_f \alpha^2 \right) \frac{1}{\alpha^2} \frac{2}{a} \sin(\omega t)  \\
        (2\lambda-{\omega}_h) \tilde{\eta}_f+{\omega}_h \frac{1}{\alpha^2} \nu(z) \\
        -{\lambda} \alpha
    \end{bmatrix}, \label{fedtautran}
\end{align}
where the function $\nu$ is defined as
\begin{align}
    \nu(z)=P(d^*+z)-P(d^*) \label{nudefined}
\end{align}
with 
\begin{align}
z = \tilde{d}_f \alpha + a \sin(\omega t)\alpha,
\end{align}
based on the identity $d = \hat{d} + a \sin(\omega t) \alpha$ and transformation \eqref{transf1}.

From Assumption \ref{assconvex}, we have
\begin{align}
    \nu(0)=0, \quad \frac{\partial }{\partial z}\nu(0)=0, \quad \frac{\partial^2 }{\partial z^2}\nu(0)=h<0. \label{nupartials}
\end{align}

\textbf{Step 2: Verification of the feasibility of \eqref{fedtautran} for averaging.}
We express the system \eqref{fedtautran} in the fast time scale $\tau = \omega t$ as follows:
\begin{align}
    &\frac{d}{d \tau}\begin{bmatrix}
    \tilde{d}_f & \hat{g}_f & \tilde{\eta}_f & \alpha
    \end{bmatrix}^T \nonumber \\
    &=\frac{1}{\omega} \begin{bmatrix}
        {\lambda} \tilde{d}_f+{k} \hat{g}_f \\ ({\lambda}-{\omega}_l) \hat{g}_f +{\omega}_l \left(\nu(z)-\tilde{\eta}_f \alpha^2 \right)  \frac{1}{\alpha^2} \frac{2}{a} \sin(\tau)   \\
        (2{\lambda}-{\omega}_h) \tilde{\eta}_f +{\omega}_h \frac{1}{\alpha^2} \nu(z) \\
        -{\lambda} \alpha
    \end{bmatrix}.  \label{systautrans}
\end{align}
We now rewrite \eqref{systautrans} compactly as:
\begin{align}
    \frac{d{\zeta}_f}{d \tau}={}(1/\omega) \mathcal{F}(\tau,\zeta_f),
\end{align}
where 
\begin{align}
    \zeta_f = \begin{bmatrix} \tilde{d}_f & \hat{g}_f & \tilde{\eta}_f & \alpha \end{bmatrix}^T,
\end{align}
and $\mathcal{F}(\tau, \zeta_f)$ denotes the vector-valued function on the right-hand side of \eqref{systautrans}.

As rigorously established in \cite{YILMAZ2025112392}, the function $\mathcal{F}(\tau, \zeta_f)$ is continuous and bounded in both arguments. Therefore, it satisfies the regularity conditions required by the averaging theorem in \cite{khalil}.

\textbf{Step 3: Averaging operation.} 
The averaged dynamics of the system \eqref{systautrans} over the period $\Pi=2\pi/\omega$ are given by:
\begin{align}
   & \frac{d}{d \tau}\begin{bmatrix}
    \tilde{d}_f^a \\ \hat{g}_f^a  \\ \tilde{\eta}_f^a 
    \\ \alpha^a
    \end{bmatrix}={}\frac{1}{\omega} \begin{bmatrix}
        {\lambda} \tilde{d}_f^a+{k} \hat{g}_f^a \\ ({\lambda}-{\omega}_l) \hat{g}_f^a \\
        (2{\lambda}-{\omega}_h) \tilde{\eta}_f^a  \\
        -{\lambda} \alpha^a
    \end{bmatrix}  \nonumber \\
    &\qquad +  \frac{1}{\omega}  \begin{bmatrix} 0 \\ {\omega}_l \frac{1}{\Pi} \int_{0}^{\Pi} \nu(\tilde{d}_f^a \alpha^a+a\sin(\sigma)\alpha^a)  \frac{1}{(\alpha^a)^2} \frac{2}{a} \sin(\sigma) d\sigma \\ {\omega}_h \frac{1}{\Pi} \int_{0}^{\Pi}  \nu(\tilde{d}_f^a \alpha^a+{a} \sin(\sigma)\alpha^a) \frac{1}{(\alpha^a)^2} d\sigma \\ 0 \end{bmatrix}, \label{averagesys}
\end{align}
where $\tilde{d}_f^a$, $\hat{g}_f^a$, $\tilde{\eta}_f^a$, and $\alpha^a$ denote the averaged versions of the corresponding state variables.

The equilibrium of the averaged system \eqref{averagesys} is then obtained as:
\begin{align}
&\begin{bmatrix}
    \tilde{d}_f^{a,e} & \hat{g}_f^{a,e} & \tilde{\eta}_f^{a,e} & \alpha^{a,e}
\end{bmatrix}^T \nonumber \\
    &\qquad =\begin{bmatrix}
    0_{1 \times n} & 0_{1 \times n} & \frac{{\omega}_h}{4({\omega}_h-2{\lambda})} h a^2 & 0
\end{bmatrix}^T, \label{aveq}
\end{align}
provided the following conditions hold:
\begin{align}
    {\omega_l}\neq{}{\lambda}, \quad
    {\omega}_h\neq{}2{\lambda}, \quad
    k \neq{} {\lambda}({\lambda}-{\omega}_l) / ( {{\omega}}_lh ).
\end{align}

\begin{figure*}[ht]
    \centering
    \includegraphics[width=1.4\columnwidth]{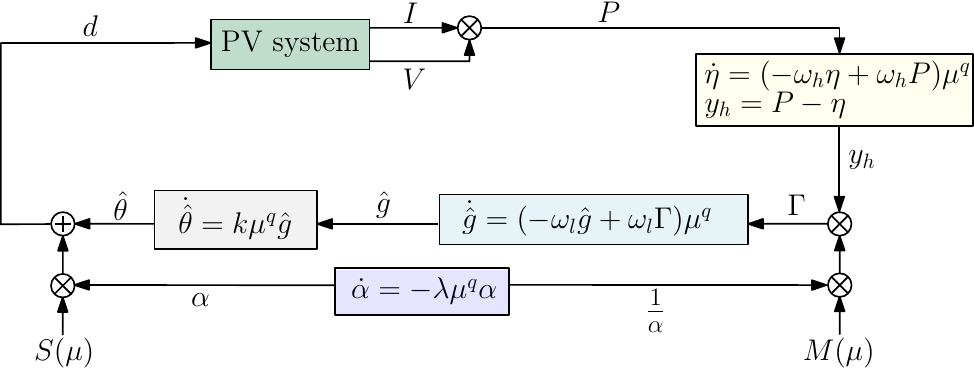}
    \caption{uPT-ES scheme. This design extends the exponential uES shown in Fig.~\ref{ESBlock} by incorporating the factor $\mu^q$, with $q \geq 1$, into all system dynamics, and by replacing the perturbation and demodulation signals with hyperbolic chirps.}
    \label{PTESBlock}
\end{figure*}

\textbf{Step 4: Stability analysis.}
The Jacobian of the averaged system \eqref{averagesys} evaluated at the equilibrium point \eqref{aveq} is given by:
\begin{align}
    &J_f^a \nonumber \\
    =&\frac{1}{\omega} \begin{bmatrix} {\lambda} & {k} & 0 & 0 \\ {\omega}_l  h & ({\lambda}-{\omega}_l) & 0 & \frac{{\omega}_l}{\Pi}  \int_0^{\Pi} \frac{\partial \left( \frac{\nu}{(\alpha^{a})^2} \frac{2}{a}\sin(\sigma)\right)}{\partial \alpha^a}  d\sigma 
    \\ 0 & 0 & (2{\lambda}-{\omega}_h) & \frac{{\omega}_h}{\Pi}  \int_0^{\Pi} \frac{\partial \left( \frac{\nu}{(\alpha^{a})^2}\right)}{\partial \alpha^a}  d\sigma  \\ 0 & 0 & 0 & -{\lambda} \end{bmatrix}.
\end{align}
Note that $J_f^a$ is block-upper-triangular and Hurwitz under the conditions given in \eqref{cond1} and \eqref{cond2}.
This establishes the local exponential stability of the averaged system \eqref{averagesys}. Applying the averaging theorem \cite{khalil}, we conclude that there exists a threshold frequency $\bar{\omega}$ such that for any $\omega > \bar{\omega}$, the original system \eqref{systautrans} admits a unique, exponentially stable periodic solution
\begin{align}
    (\tilde{d}_f^{\Pi}(\tau), \hat{g}_f^{\Pi}(\tau), \tilde{\eta}_f^{\Pi}(\tau), \alpha^{\Pi}(\tau)) 
\end{align}
with period $\Pi=2\pi/\omega$, and this solution satisfies the bound:
\begin{align}
    \left| \begin{bmatrix} \tilde{d}_f^{\Pi}(\tau) \\ \hat{g}_f^{\Pi} (\tau)  \\ \tilde{\eta}_f^{\Pi}(\tau)-\frac{{\omega}_h}{4({\omega}_h-2{\lambda})} h a^2  \\ \alpha^{\Pi}(\tau)  \end{bmatrix} \right| \leq \mathcal{O}\left( \frac{1}{\omega} \right).
\end{align}
In other words, all solutions $(\tilde{d}_f(\tau), \hat{g}_f(\tau),  \tilde{\eta}_f(\tau), \alpha(\tau))$ exponentially converge to an $\mathcal{O}\left( 1/\omega \right)$-neighborhood of the equilibrium point in \eqref{aveq}. The signal $\alpha(\tau)$, in particular, exponentially converges to zero. Recalling the transformations \eqref{trans}, we can deduce that the system \eqref{ffed} with states $\tilde{d}(t), \hat{g}(t), \tilde{\eta}(t)$ has a unique solution and is exponentially stable at the origin. Specifically:
\begin{itemize}
\item $\tilde{d}(t)$ and $\hat{g}(t)$ decay exponentially with rate $\lambda$,
\item $\tilde{\eta}(t)$ decays exponentially with rate $2\lambda$.
\end{itemize}

\textbf{Step 5: Convergence to extremum.} 
Building on the results from Step 4 and recalling from \eqref{trans} and Fig.~\ref{ESBlock} that
\begin{align}
    d(t)={}&\alpha(t) \tilde{d}_f(t)+d^*+\alpha(t) a \sin(\omega t),
\end{align}
we conclude that $d(t)$  converges exponentially to the optimal point $d^*$ at rate $\lambda$, since $\alpha(t) \to 0$ exponentially and $\tilde{d}_f(t)$ remains bounded for all $t \geq t_0$. Moreover, using the definition of $\nu$ in \eqref{nudefined} and \eqref{nupartials}, it follows that the output signal $P(d(t))$  converges exponentially to the optimal value $P(d^*)$ with rate $2\lambda$. This concludes the proof of Theorem~\ref{the1}.
\hfill \end{proof}

\section{Unbiased Prescribed-Time Extremum Seeker} \label{PTESSec}

To further accelerate convergence to the maximum power point, we propose a novel ES algorithm termed unbiased prescribed-time extremum seeking (uPT-ES). As its name suggests, uPT-ES ensures unbiased convergence to the optimum within a user-defined terminal time $T$. The structure of uPT-ES is illustrated schematically in Fig. \ref{PTESBlock}. In contrast to the exponentially convergent $\alpha$-dynamics used in \eqref{mu}, we define a prescribed-time convergent $\alpha$-dynamics as 
\begin{align}
    \dot{\alpha}(t)=&{}-\lambda \mu^q(t-t_0) \alpha(t), \qquad \alpha(t_0)=\alpha_0 \label{mudot}
\end{align}
with $q \geq 1$ and the time-varying scaling function $\mu(\cdot)$ is given by:
\begin{align}
    \mu(t-t_0)={}\frac{T}{T+t_0-t}, \quad t \in [t_0,t_0+T). \label{xidef}
\end{align}
The solution to~\eqref{mudot} is:
\begin{align}
    \alpha(t)={}&\alpha(t_0) e^{-\lambda\int_{t_0}^t \mu^q(\sigma-t_0)d\sigma} \nonumber \\
    ={}&\begin{cases}
    \alpha_0 \mu^{-\lambda T}(t-t_0),& \text{if } q= 1,\\
    \alpha_0 e^{- \frac{\lambda T}{q-1} \left(\mu^{q-1}(t-t_0)-1\right)},  & \text{if } q> 1 
\end{cases} \label{mudef}
\end{align}
with the key property that $\alpha(t_0+T)=0$. Increasing the value of $q$ leads to a more rapid decay of $\alpha(t)$, due to the steeper growth of $\mu(t-t_0)$.

The closed-loop dynamics of the uPT-ES system, shown in Fig.~\ref{PTESBlock}, can be written as:
\begin{align}
    \frac{d}{d t}\begin{bmatrix}
    \tilde{d} \\ \hat{g} \\ \tilde{\eta} \\ \alpha
    \end{bmatrix}=\begin{bmatrix}
        k \mu^q \hat{g} \\ -\omega_l \mu^q \hat{g}+\omega_l \mu^q (P-P(d^*)-\tilde{\eta}) \frac{1}{\alpha} M(\mu) \\
        -\omega_h \mu^q \tilde{\eta}+\omega_h \mu^q (P-P(d^*)) \\
        -\lambda \mu^q \alpha
    \end{bmatrix}, \label{ptclosed}
\end{align}
in view of the transformations \eqref{transf1}, \eqref{transf2}, and $\eta$ is governed by
\begin{align}
    \dot{\eta}=(-\omega_h  \eta+\omega_h P)\mu^q.
\end{align}

To ensure prescribed-time convergence to the extremum, we replace the conventional sinusoidal perturbations with time-varying (or ``chirpy'') signals whose frequency increases with time:
\begin{align}
    S(\mu)
    &={}\begin{cases}
    a \sin(\omega (t_0+T\ln(\mu)) , & \text{if } q= 1, \\
    a \sin\left(\omega  \left(t_0+T \frac{\mu^{q-1}-1}{q-1}  \right)  \right), & \text{if } q> 1, \label{Smudef}
\end{cases}  \\
    M(\mu) 
    &={}\begin{cases}
    \frac{2}{a} \sin(\omega (t_0+T\ln(\mu)) , & \text{if } q= 1, \\
    \frac{2}{a} \sin\left(\omega  \left(t_0+T \frac{\mu^{q-1}-1}{q-1}  \right)  \right), & \text{if } q> 1. \label{Mmudef}
\end{cases}
\end{align}

\begin{figure}[t]
    \centering
    \includegraphics[width=\columnwidth]{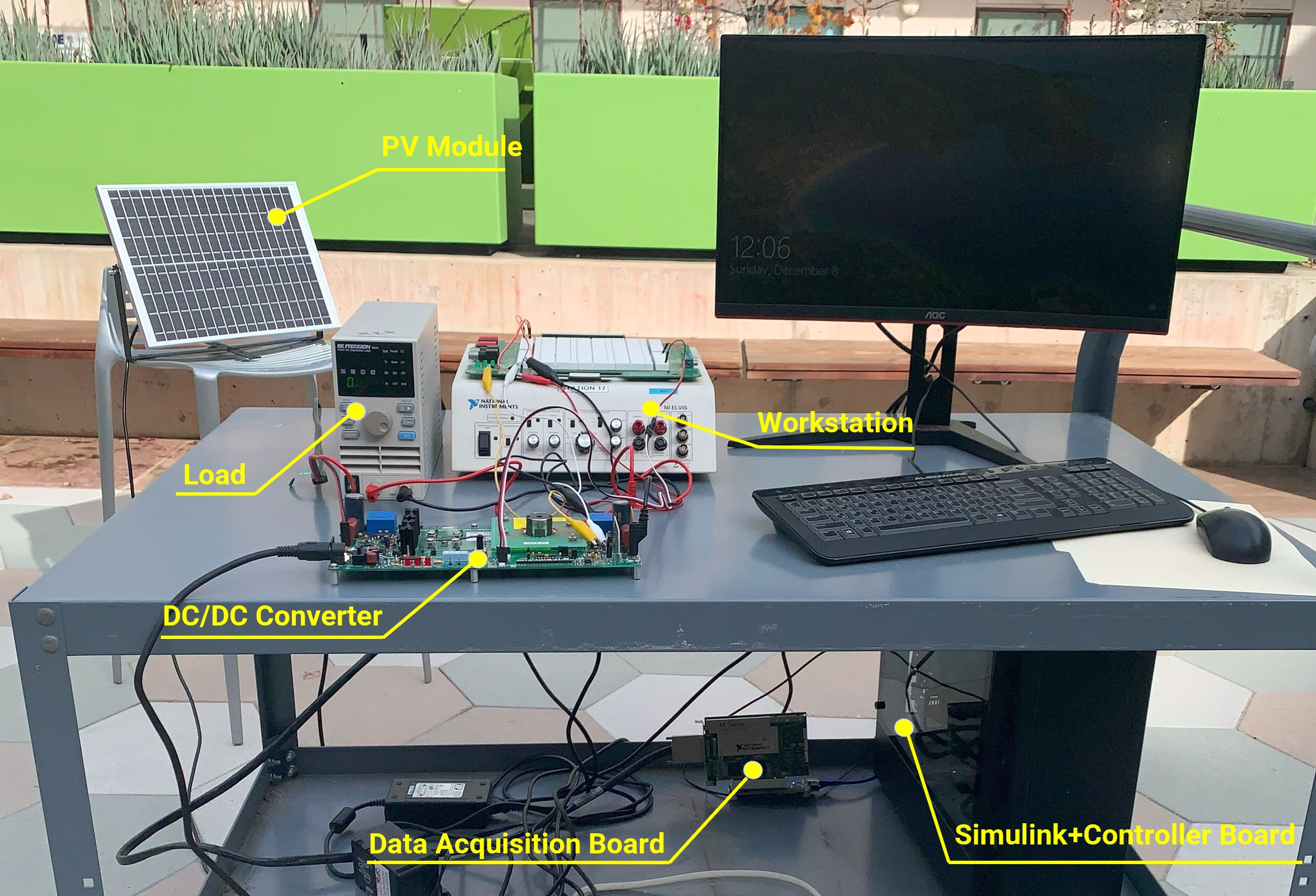}
    \caption{The experimental setup of the PV system.}
    \label{setup}
\end{figure}

\begin{figure*}[ht]
    \centering
    \begin{subfigure}[t]{0.32\textwidth}
        \centering
        \includegraphics[width=\linewidth]{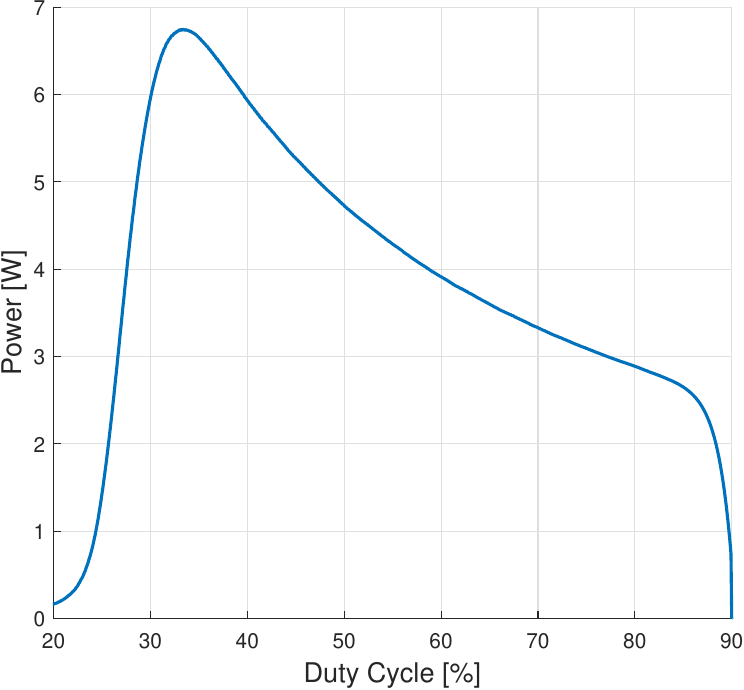}
        \caption{}
        \label{uesfiga}
    \end{subfigure}
    \hfill
    \begin{subfigure}[t]{0.32\textwidth}
        \centering
        \includegraphics[width=\linewidth]{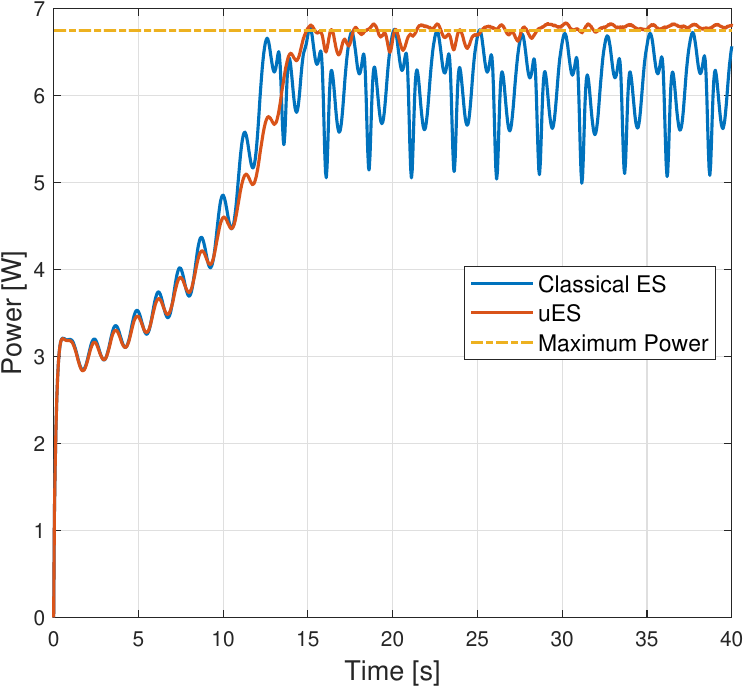}
        \caption{}
        \label{uesfigb}
    \end{subfigure}
    \hfill
    \begin{subfigure}[t]{0.32\textwidth}
        \centering
        \includegraphics[width=\linewidth]{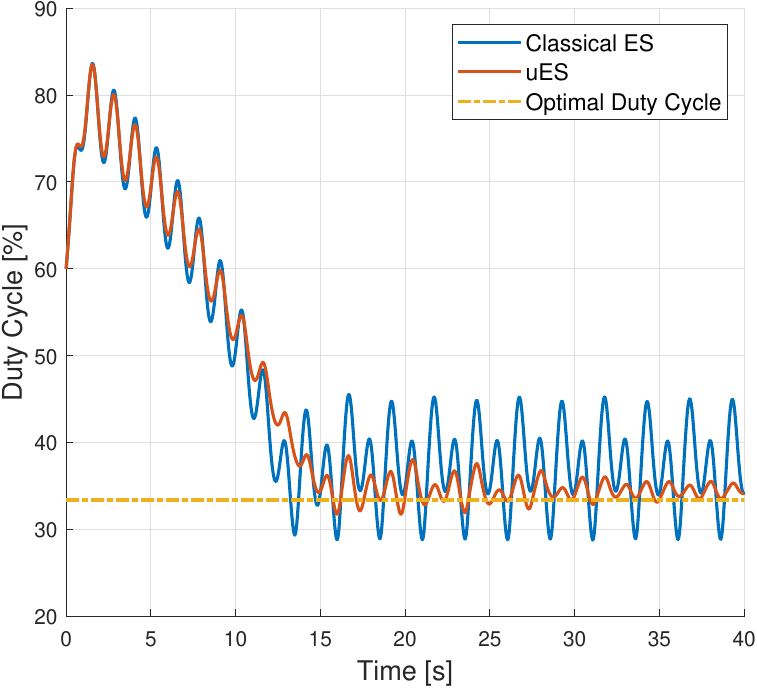}
        \caption{}
        \label{uesfigc}
    \end{subfigure}
    \caption{Experimental validation of the proposed uES for MPPT, compared to classical ES: (a) Power vs. duty cycle; (b) Power convergence; (c) Duty cycle convergence.}
    \label{fig:uES_validation}
\end{figure*}

The following theorem is reproduced from \cite{YILMAZ2025112392}.
\begin{theorem} \label{the2}
Consider the feedback system \eqref{ptclosed} with the parameters that satisfy \eqref{cond1}, \eqref{cond2}. There exists 
$\bar{\omega}$ and for any $\omega > \bar{\omega}$ 
there exists an open ball $\mathcal{B}$ centered at the point $(\hat{d}, \hat{g}, {\eta}, \alpha)=(d^*, 0, P(d^*),0) = : \Upsilon$ such that  for any initial condition starting in the ball $\mathcal{B}$, the system \eqref{ptclosed} has a unique solution and the solution converges to $\Upsilon$ in prescribed time $T$. Furthermore, $P(d(t))$ converges to $P(d^*)$ in prescribed time $T$.
\end{theorem}

A brief sketch of the proof is provided below:

\begin{proof}

\textbf{Step 0: Time dilation from $t$ to $\check{\tau}$.}
Define the following time dilation and contraction transformations:
\begin{align}
    \check{\tau}={}&\begin{cases}
    t_0+T\ln(\mu), & \text{if } q= 1,\\
    t_0+T\left( \frac{\mu^{q-1}-1}{q-1}  \right),  & \text{if } q> 1, 
\end{cases}  \label{dilat}  \\
    t={}&\begin{cases}
    t_0+T\left(1-e^{-\frac{\check{\tau}-t_0}{T}}\right), & \text{if } q= 1,\\
    t_0+T\left(1-\left(\frac{T}{T+(q-1)(\check{\tau}-t_0)} \right)^{\frac{1}{q-1}} \right), & \text{if } q> 1, 
\end{cases} \label{contrac}
\end{align}
for $\check{\tau} \in [t_0,\infty)$, $t \in [t_0,t_0+T)$. 
Differentiating \eqref{dilat} with respect to $t$, we obtain
\begin{align}
    \frac{d \check{\tau}}{dt}={}\mu^q(t-t_0), \label{dtaudt}
\end{align}

Using this time dilation, we rewrite the system \eqref{ptclosed} in the $\check{\tau}$-domain as:
\begin{align}
    \frac{d}{d \check{\tau}}\begin{bmatrix}
    \tilde{d} \\ \hat{g} \\ \tilde{\eta} \\ \alpha
    \end{bmatrix}=\begin{bmatrix}
        k  \hat{g} \\ -\omega_l \hat{g}+\omega_l  (P-P(d^*)-\tilde{\eta}) \frac{1}{\alpha}\mathcal{M}(\check{\tau})  \\
        -\omega_h \tilde{\eta}+\omega_h (P-P(d^*)) \\
        -\lambda \alpha
    \end{bmatrix}, \label{ptclosedtaucheck}
\end{align}
with the perturbation and demodulation signals defined as
\begin{align}
    \mathcal{S}(\check{\tau})=&  a \sin(\omega \check{\tau}), \label{Sttau} \\
    \mathcal{M}(\check{\tau})=&\frac{2}{a} \sin(\omega \check{\tau}). \label{Mttau}
\end{align}

Note that the system \eqref{ptclosedtaucheck}, together with \eqref{Sttau} and \eqref{Mttau}, has the same structure as the system \eqref{ffed}, except that it evolves with respect to the dilated time variable $\check{\tau}$  instead of the original time $t$. This time transformation plays a critical role by making the perturbation and demodulation signals periodic with constant frequencies in the $\check{\tau}$-domain, and this allows the direct application of the averaging theorem. The remainder of the proof then follows by applying Steps 1 through 5 from the proof of Theorem~\ref{the1}, in the $\check{\tau}$-domain. Finally,  applying the transformation \eqref{contrac} contracts time back from $\check{\tau}$ to $t$ and completes the proof.
\hfill
\end{proof}

\begin{figure*}[t]
    \centering
    \includegraphics[width=2\columnwidth]{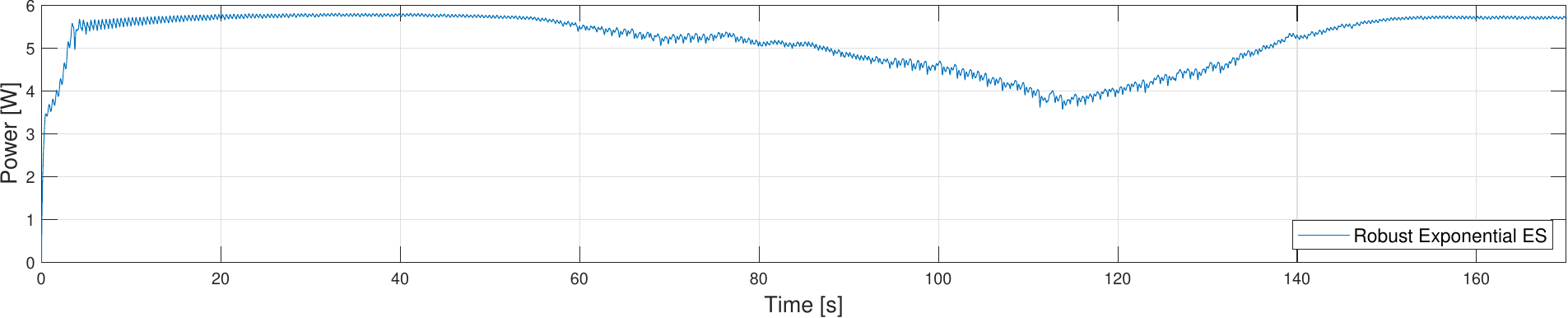}
    \caption{Tracking the shifting optimum of a PV system under shading using uES.}
    \label{result3}
\end{figure*}

\section{Experimental Setup and Results}
The experimental validation of the proposed control algorithm is conducted on a hardware testbed. The setup, shown in Fig. \ref{setup}, consists of the certain components. A commercially available PV module is mounted outdoors to receive direct sunlight. The module serves as the power source and exhibits nonlinear power-voltage characteristics dependent on irradiance and temperature. A DC/DC buck-boost converter is used to regulate the output voltage of the PV module. The converter interface allows for real-time perturbations of duty cycle. The load simulates a power demand and operates in constant voltage. National Instruments Data Acquisition Board (NI-DAQ) hardware is used to collect measurements of voltage, current, and power. These signals are digitized and fed to the computer for real-time analysis.  The controller is implemented on a real-time embedded board programmed using MATLAB/Simulink. This allows deployment of the control algorithm from simulation to hardware and supports closed-loop operation. The host computer runs the Simulink model and interfaces with both the controller board and the DAQ system. It enables logging of experimental data and adjustment of parameters during testing.

Due to sensor noise and potential shifts in the maximum power point and optimal duty cycle under varying irradiance, implementing unbiased designs presents practical challenges. To address this, we modify the oscillation amplitude dynamics during implementation as follows:
\begin{align}
\dot{\alpha}(t) = -\lambda \alpha(t) + \lambda \beta, \qquad \alpha(t_0) = \alpha_0, \label{mumodif}
\end{align}
which prevents $\alpha(t)$ from decaying to zero and instead drives it to a small value $\beta$. This compromises strict unbiasedness but improves robustness in real-world conditions.

The experimental results in Fig. \ref{fig:uES_validation} demonstrate the effectiveness of the proposed uES algorithm, with parameters set to $k=0.01, \lambda=0.05, a=0.2, \omega=5, \omega_h=\omega_l=3, \alpha_0=1, \beta= 0.1$. The static power curve, shown in Fig. \ref{uesfiga}, exhibits a clear unimodal behavior with a peak around a 34\% duty cycle, confirming the presence of a well-defined maximum power point. Fig. \ref{uesfigb} presents the measured power over time using both the classical ES and the proposed uES schemes. While both methods ultimately reach the vicinity of the maximum, uES achieves convergence with smaller amplitude oscillations around the optimum. Fig. \ref{uesfigc} shows the evolution of the duty cycle during the experiment. The uES controller drives the duty cycle toward the optimal value with minimal oscillation, in contrast to the classical ES which exhibits much larger dithering.

As discussed in Fig. \ref{fig:exponent_uES_track}, changes in irradiance cause the optimal value to shift over time. The ability of the developed design for tracking the shifting optimum is tested, and results are provided in Fig. \ref{result3}. The PV panel is started shading at $t=40$ seconds, gradually shaded until around $t=115$ seconds, and shade is uncovered gradually for the remainder of the test. The developed design achieve a close tracking of the shifting optimum, adjusting the duty cycle.

\begin{figure*}[t]
    \centering
    \begin{subfigure}[t]{0.3\textwidth}
        \centering
        \includegraphics[width=\linewidth]{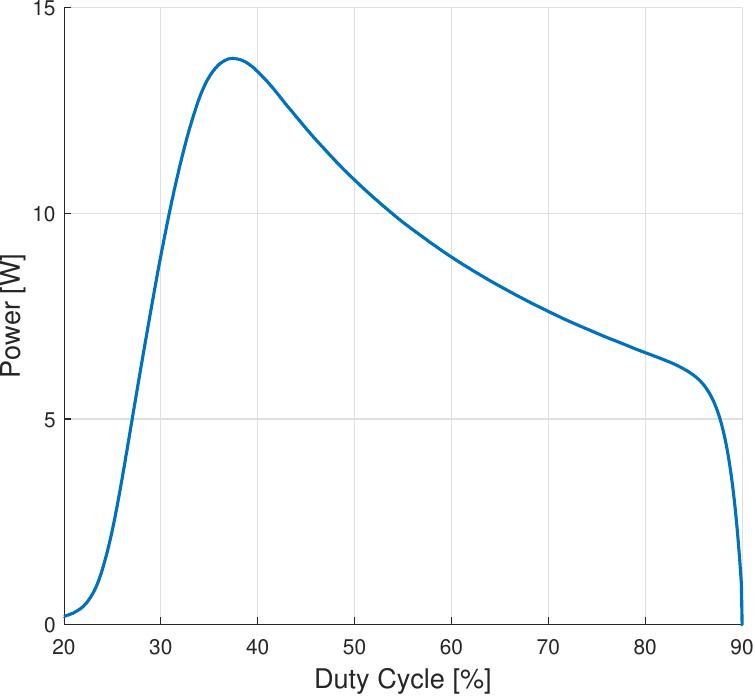}
        \caption{}
        \label{uptesfiga}
    \end{subfigure}
    \hfill
    \begin{subfigure}[t]{0.3\textwidth}
        \centering
        \includegraphics[width=\linewidth]{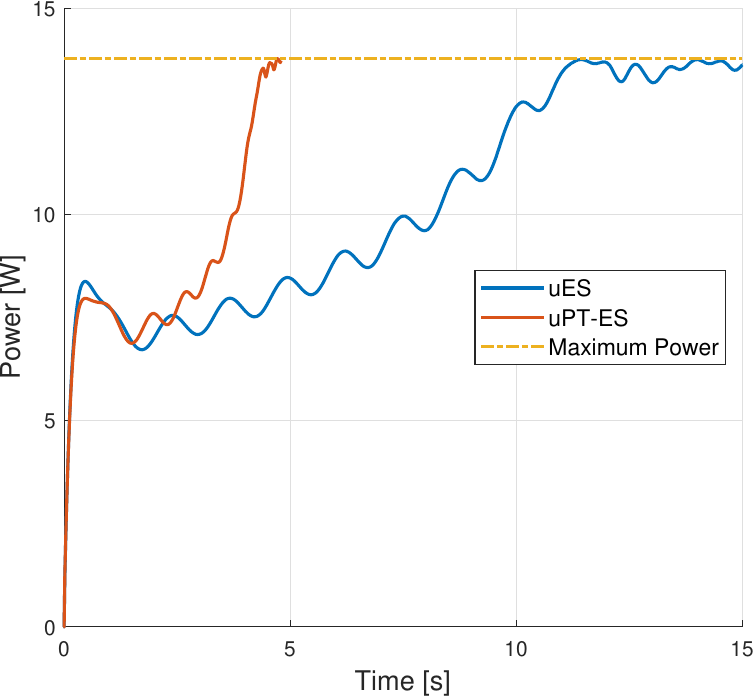}
        \caption{}
        \label{uptesfigb}
    \end{subfigure}
    \hfill
    \begin{subfigure}[t]{0.3\textwidth}
        \centering
        \includegraphics[width=\linewidth]{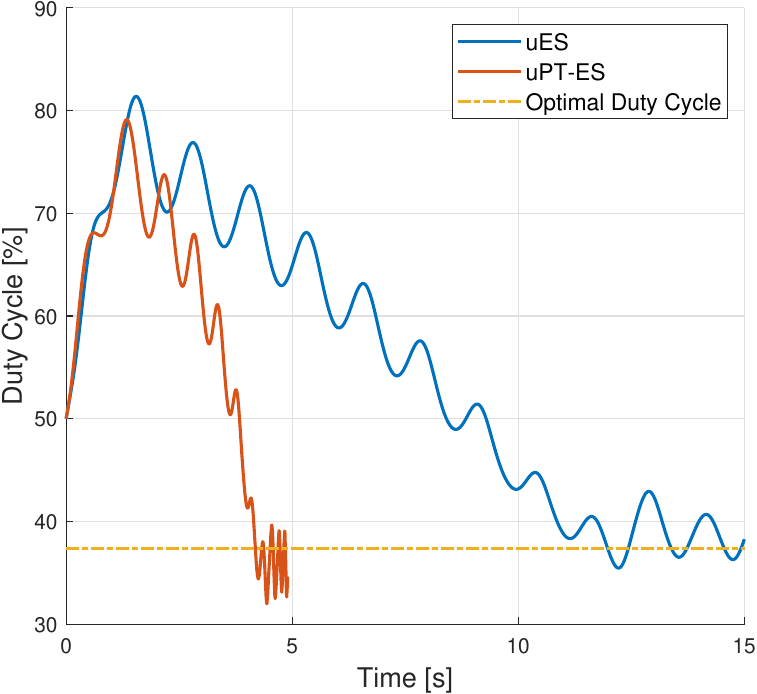}
        \caption{}
        \label{uptesfigc}
    \end{subfigure}
    \caption{Experimental comparison of uPT-ES and uES under a user-prescribed terminal time of 5 seconds: (a) Power vs. duty cycle; (b) Power convergence; (c) Duty cycle convergence.}
    \label{result2}
\end{figure*}

The results in Fig. \ref{result2} demonstrate the performance advantage of the proposed uPT-ES algorithm over the unbiased uES. The terminal time is set to 6 seconds, while the experiment is intentionally stopped at $t=5$ seconds to prevent excessive growth of the blow-up function $\mu(t)$. The parameters are chosen as 
$q=1, T=6$ s, $\omega=5$ rad/s. Fig. \ref{uptesfiga} presents the static power curve, and Fig. \ref{uptesfigb} illustrates the accelerated power convergence achieved with uPT-ES. This acceleration is enabled by the use of a chirp signal in the perturbation, as clearly reflected in the duty cycle trajectory shown in Fig. \ref{uptesfigc}.

\section{Conclusion}
This paper has proposed and experimentally validated two novel model-free MPPT algorithms: exponential unbiased extremum seeking and unbiased prescribed-time extremum seeking. The algorithms overcome the limitations of classical ES methods. The uES design eliminates steady-state bias while preserving exponential convergence, and the uPT-ES framework further accelerates convergence by guaranteeing finite-time settling to the MPP. Theoretical convergence guarantees for both static and dynamic optima are established and supported by extensive experimental results. In practical implementations, minor modifications to the perturbation amplitude dynamics are shown to enhance robustness under measurement noise and environmental variability. Extensions to distributed PV systems remain promising directions for future research.

\section*{References}
\vspace{-0.5cm}

\end{document}